\documentclass[twoside]{amsart}

\usepackage[latin1]{inputenc}
\usepackage{bbold}
\usepackage{float}
\usepackage{url}

\newtheorem{thm}{Theorem}[section]

\newtheorem{question}[thm]{Question}

\floatstyle{ruled}
\newfloat{example}{tbph}{loe}[section]
\floatname{example}{Example}

\numberwithin{equation}{section}

 \newcommand{\Cpp}%
  {C\nolinebreak\hspace{-.05em}\raisebox{.4ex}{\tiny\bf +}%
    \nolinebreak\hspace{-.10em}\raisebox{.4ex}{\tiny\bf +}}
                                                                               
\newcommand\dterm[1]{{\em #1}}
\newcommand\ideal[1]{\langle #1 \rangle}
\newcommand\lcm{\operatorname{lcm}}
\newcommand\meet\wedge
\newcommand\join\vee
\newcommand\bigmeet\bigwedge
\newcommand\bigjoin\bigvee
\newcommand\iso\cong
\newcommand\Tor{\operatorname{Tor}}

\newcommand\ZZ{\mathbb Z}
\newcommand\QQ{\mathbb Q}
\newcommand\NN{\mathbb N}

\newcommand\ie{i.e.\@ }

\begin{document}

\author{Mikael Johansson}
\thanks{Submitted to Rendiconti del Seminario Matematico
  dell'Universita` e del Politecnico di Torino, for publication in the
  proceedings of the School and Workshop on Algebraic Geometry and
  Statistics}
\title[Computation of Poincaré-Betti series]{Computation of Poincaré-Betti series for monomial rings}
\address{Department of Mathematics, Stockholm University, SE-10691 Stockholm}
\email{mik@math.su.se}

\begin{abstract}
 
The multigraded Poincaré-Betti series $P_R^k(\mathbb x; t)$ of a monomial ring $k[\mathbb x]/\ideal M$ on a finite number of monomial
generators has the form $\prod_{x_i\in\mathbb
x}(1+x_it)/b_{R,k}(\mathbb x;t)$, where $b_{R,k}(\mathbb x;t)$ is a
polynomial depending only on the monomial set $M$ and the
characteristic of the field $k$. I present a computer program designed
to calculate the polynomial $b_{R,k}$ for a given field characteristic
and a given set of monomial generators.
\end{abstract}

\maketitle

\section{Introduction}

Let $Q=k[\mathbb x]=k[x_1,\dots,x_r]$ be the polynomial ring over a
field $k$ with $r$ variables. The ring has a natural
$\NN^r$-grading\footnote{Throughout this paper, $0\in\NN$.}
by setting $\deg(x_i)=e_i$ for the canonical basis vectors $e_i$ of
$\NN^r$. Setting $|(a_1,\dots,a_r)|=a_1+\dots+a_r$, we can derive an
$\NN$-grading of $Q$ from this $\NN^r$-grading. We write $\deg(m)$ for
the $\NN^r$-degree of a monomial $m$, and given
$(a_1,\dots,a_r)=\alpha\in\NN^r$, we write
$x^\alpha=x_1^{a_1}x_2^{a_2}\dots x_r^{a_r}$

The $\NN^r$-grading and its inherent $\NN$-grading both are inherited
from $Q$ to the ring $R=Q/I$ where $I$ is a \dterm{monomial ideal},
\ie an ideal generated by monomials in $Q$. We call $R$ a
\dterm{monomial ring}.

The $\NN^r$-grading inherits, via minimal resolutions respecting to
the grading, to an $\NN^r$-grading on $\Tor^R(k,k)$. Thus, we can define
the multigraded Poincaré-Betti series $P^R_k(\mathbb x; t)$ of an
$\NN^r$-graded ring $R$ over $k$: 
$$\sum_{i\in\NN}\sum_{\alpha\in\NN^r}\dim_k\Tor^R_i(k,k)_\alpha
x^\alpha t^i$$

From $P_k^R(\mathbb x;t)$, the simple Poincaré-Betti series can
be calculated as $P_k^R(t)=P_k^R(1,\dots,1;t)$.

Since \cite{joeb-rational}, it is known that the Poincaré-Betti series
of a monomial ring $R$ holds the form 
$P^R_k(t)=\frac{(1+t)^n}{b_{R,k}(t)}$ for some polynomial
$b_{R,k}(t)$.

Moreover, Alexander Berglund proved in \cite{alex-arxiv} that $\deg
b_{R,k}(t)<2n$ for monomial rings with $n$ monomial generators. It
follows that there are finitely many Poincaré-Betti serie occurring
for a fixed $n$ at all.
 
We may define a partially ordered graph, or po-graph, to be a graph
with a partial order on the vertices. Two po-graphs are said to be
isomorphic if there is a simultaneous isomorphism of the graph and the
partial order.

Let $Q=k[\mathbb x]$ and $Q'=k[\mathbb x']$ be polynomial rings for
two finite variable sets $\mathbb x=\{x_1,\dots,x_r\}$ and $\mathbb
x'=\{x_1,\dots,x_{r'}\}$, and let $I$ and $I'$ be monomial ideals in
$Q$ and $Q'$ respectively, with $M$ and $M'$ the sets of generators
for each ideal. For some set $S$ of monomials, we denote by
$L_S$ the set of all least common multiples of subsets of $S$. $L_S$
can be equipped with the structure of a po-graph, ordering monomials
by divisibility and adding an edge between two elements when they have
a non-trivial common factor.

Luchezar Avramov shows in \cite{avramov} that if two rings $R=Q/I$ and
$R=Q'/I'$, with $I$ and $I'$ generated by the monomial sets $M$ and
$M'$ respectively are such that $L_M\iso L_{M'}$, then
$b_{R,k}(t)=b_{R',k}(t)$. From this follows that for a fixed field and
a fixed number of monomial generators, only finitely many different
Poincaré-Betti series can occur. Avramov further proves that the
limitation to a fixed field is superfluous.

In \cite{alex-arxiv}, Alexander Berglund, proving the conjecture by
Avramov that $\deg(b_{R,k}(t))\leq 2n$ whenever the monomial ideal
generating the ring $R$ has $n$ generator, constructs combinatorially
a minimal model for $R$ and gives a characterisation of the
Poincaré-Betti series denominator polynomial in terms of the homology
of associated simplicial complexes. Avramov's earlier observation that
only finitely many different Poincaré-Betti series exist for a fixed
number of generators for the monomial ideal $I$ follows as an
immediate consequence of Berglund's construction.

I have in the course of my M.Sc.\@ thesis work \cite{mik-master} continued
Berglund's work by implementing his formula in a program, written in
\Cpp, capable of calculating simplicial
homology over fields of arbritrary characteristic, as well as
explicitly calculating $b_{R,k}(t)$ for arbritrary characteristic of
the coefficient field $k$.

\subsection{Simplicial Complexes}

Since we will work a lot with simplicial homology, I will take a few
moments to review definitions and terminology. A simplicial complex on
a set $V$ is a set $\Delta$ of subsets of $V$ such that if
$G\in\Delta$ and $F\subset G$, then $F\in\Delta$. $V$ is called the
\dterm{vertex set} of $\Delta$. All simplicial complexes I shall refer
to will have $V=\bigcup\Delta$ unless otherwise stated. The
$i$-\dterm{faces} or $i$-\dterm{simplices} of $\Delta$ are precisely
the elements in $\Delta$ of cardinality $i+1$. 

To a simplicial complex $\Delta$ we can associate an augmented chain
complex $\widetilde C(\Delta)$ with $\widetilde C_i(\Delta)$ the free
abelian group on the $i$-faces of $\Delta$, where we consider
$\emptyset$ to be the unique $-1$-simplex. We equip $\widetilde
C(\Delta)$ with the standard differential of degree $-1$. Thus
$$H_i(\widetilde C(\Delta))=\widetilde H_i(\Delta)$$

As usual, for an abelian group $G$, we set $\widetilde C(\Delta;
G)=\widetilde C(\Delta)\otimes_{\ZZ}G$ and $\widetilde
H_i(\Delta;G)=H_i(\widetilde C(\Delta;G))$. 

For a simplicial complex $\Delta$ we define the Alexander dual
$$\Delta^\join=\{F\subseteq V\mid V\setminus F\not\in\Delta\}$$

For simplicial complexes $\Delta$ and $\Delta'$ with disjoint vertex
sets we define the join
$$\Delta*\Delta'=\{F\cup F'\mid F\in\Delta, F'\in\Delta'\}$$
and the dual join
$$\Delta\cdot\Delta'=(\Delta^\join*\Delta^{'\join})^\join$$
and note that then
$$\Delta\cdot\Delta'=\{F\in V\cup V'\mid F\cap V\in\Delta\mbox{ or }
F\cap V'\in\Delta'\}$$

In \cite[Lemma 5.5.3]{cohen-macaulay}, it's shown that 
$$\widetilde H_i(\Delta; k)\iso\widetilde H_{n-3-i}(\Delta^\join;k)$$
for a complex $\Delta$ on $n$ vertices.

For a graded vector space $H=\bigoplus_{i\in\NN}H_i$, we will write
$H(t)$ for the \dterm{generating function} $\sum_{i\in\NN}\dim_k
H_it^i$ of $H$. We then can find that the join of complexes as well as
the Alexander dual gives rise to rather easily handled equalities on
the level of generating functions for their respective homologies,
again with $\Delta$ a complex on $n$ vertices, and $\Delta'$ some
other complex.
\begin{align}
t^n\widetilde H(\Delta^\join;k)(t^{-1}) 
  &= t^3\widetilde H(\Delta;k)(t) \\
\widetilde H(\Delta*\Delta';k)(t) 
  &= t\widetilde H(\Delta;k)(t)\cdot\widetilde H(\Delta';k)(t)
\end{align}
where the $t$ factor in the latter equation comes from the fact that
a simplex with $d$ elements is considered to have dimension $d-1$.

\section{Berglund's work}

In \cite{alex-arxiv}, Berglund treats the theoretical
aspects of computation of the Poincaré-Betti series denominator for
monomial rings. I will not repeat all of his proof here, but rather
reference his work to establish the vocabulary and touch the results I
will need for my own work.
We can find
$Q=k[\mathbb x]$ is $\NN^r$-graded by assigning to a monomial
$x_1^{a_1}\dots x_r^{a_r}$ the element $\alpha=(a_1,\dots,a_r)$. We
write $x^\alpha$ for $x_1^{a_1}\dots x_r^{a_r}$. The monomial
$x^\alpha$ is said to be squarefree if $\alpha\in\{0,1\}^r$.

By a construction by Weyman and Fröberg \cite{weyman,ralff}, it is
enough to
treat squarefree monomial sets, since an easy procedure can be used to
go from a monomial ring to a squarefree monomial ring with the same
homological properties. So we can assume that $R=Q/I$ is a squarefree
monomial ring with $I$ generated by the monomial antichain $M$ of
cardinality $n$.

For a set $S$ of monomials, set $m_S=\lcm(m:m\in S)$. In particular
$m_\emptyset=1$. Thus
$L_M=\{m_S\mid S\subset M\}$. For a specific monomial $m$ and a
monomial set $M$ set $M_m=\{m'\in M\mid m'|m\}$.

Now, for $I=\ideal M$ a monomial ideal in $Q$ generated minimally by
the antichain $M$, Berglund introduces the complex
$$\Delta_M=\{S\subseteq M\mid m_S\neq m_M\mbox{ or }S\mbox{ disconnected}\}$$
where connectivity is for $S$ as a subgraph of $L_M$.

Using multigraded ring-deviations
$$\epsilon_{i,\alpha}=\dim_k\widetilde H_{i-3}(\Delta_{M_{x^\alpha}};k)$$ 
Berglund gives the squarefree part of the multigraded Poincaré-Betti
polynomial 
\begin{equation}\label{pbpoly}
b_{R,k}(\mathbb x,t)\equiv\prod_{x\alpha\in L_M}(1-x^\alpha p_\alpha(t))\pmod{\ideal{x_1^2,\dots,x_r^2}}
\end{equation}
with $p_\alpha(t)=t^3\widetilde
H(\Delta_{M_{x^\alpha}};k)(t)$. Backelin demonstrated already in
\cite{joeb-rational} that the denominator polynomial will be
squarefree whenever the monomials generating the ideal all are.

Berglund then goes on to find several more theoretically pliable forms
of this particular formula; expanding the product and taking the
irrelevance of non-squarefree terms into account, he arrives at the
form 
\begin{equation}\label{pbpoly2}
b_{R,k}(\mathbb x;t)=1+\sum_S m_S(-t)^{c(S)+2}\widetilde H(\Delta_S;k)(t)
\end{equation}
where the sum is taken over all non-empty \dterm{saturated} subsets of
$M$ and $c(S)$ counts the number of graph components of $S$ as a
subgraph to $L_M$. We define the \dterm{saturation} of a subset
$S\subset M$ as the set of all monomials in $M$ that divide the least
common multiple of some connected component of $S$ as subgraph to
$L_M$. A set is saturated if it is equal to its saturation.

\section{The resulting application -- {\tt poincare}}

My own achievement is that I have constructed a computer program to
calculate simplical complex homology and mainly to calculate the
denominator polynomials of Poincaré-Betti series using Berglunds
methods. I will devote this section to a discussion of the program,
which can be fetched in its latest version under the MIT software
license from \url{http://www.math.su.se/~mik/poincare/}.

The form deemed most promising for implementation as I started was the
form given in \eqref{pbpoly} -- mainly since the formulation in terms
of saturated subsets had at that time not yet matured. Thus, I have
implemented specific \Cpp{} classes for calculating in the ring
$Q/\ideal{x_1^2,x_2^2,\dots,x_r^2}$ and let the final product forming
the polynomial take place in that particular ring. The only part
forming any kind of complexity for the straightforward implementation is
that of forming the complex $\Delta_{M_{x^\alpha}}$ and calculating
its homology over the specified characteristic.

The construction of $\Delta_{M_{x^\alpha}}$ is done with a modified
kind of breadth-first search: monomials are stored in a queue along
with an index keeping track of which of the monomials covering the
particular monomial that have already been tried. Thus, for each
monomial in the queue, all later covering monomials are tried one
after the other, and upon compliance with the two conditions -- that
the least common multiple of all generating monomials dividing the
candidate is equal to $x^\alpha$ and that those generating monomials
are connected as a graph -- the monomial is added to the queue
carrying a testing index one higher than the index that produced
it. This algorithm does yield a speed increase compared to the earlier
algorithm that simply tested all monomials for both conditions;
but still is not optimal by far.

Once the simplicial complex as such has been constructed, the
calculation of its homology commences. This is calculated degree by
degree, constructing a matrix with entries in $\{0,\pm 1\}$ and
fetching its rank from the matrix routines in the Pari library
\cite{pari}. 

This is wrapped in a text-mode user interface, using the GNU Readline
library to facilitate command history and command editing. The user
interface reads in space-separated lists of monomials as input to the
{\tt add simplex} and {\tt add monomial} commands. A monomial, to the
program, is a {\tt *}-separated list of strings of characters, where
each separated string is taken to be the name of a variable. The
variable names must avoid \verb|+-*/^,.| and whitespace, but can use any
other characters. Any string occuring in such a position will be
interpreted as a variable and added to an internal dynamic variable
pool.

The user interface wraps, among other things, around the Weyman-Fröberg
method for conversion to a squarefree monomial ring. The conversion is
done transparently, using several internal variables that are easily
converted back to the original variables before printing the answer.
The output rendered by the program is written in such a way that other
computer algebra systems should have an easy time handling it.

As an example on how the program works, I give in example
\ref{ex:poincare} a session, calculating first the simplicial homology
of the projective plane over $\QQ$ as well as over $\ZZ_2$ and then
calculating the Poincaré-Betti denominator polynomial of
$\QQ[x,y,z]/\ideal{x^2,xy,yz}$.

\begin{example}
\begin{verbatim}
Welcome to the Poincaré calculator. You can use this program to
calculate simplicial homology over prime fields and to calculate
the denominator polynomial of the Poincaré-Betti series of monomial
rings.
(c) 2004 Mikael Johansson
This program is released under the MIT License

> add simplex a*b*e a*b*f a*c*d a*c*f a*d*e
> add simplex b*c*d b*c*e b*d*f c*e*f d*e*f
> homology
Calculating homology ranks...
*****  Hilbert series of simplicial homology *****
0
> char 2
New characteristic: 2
> homology
Calculating homology ranks...
*****  Hilbert series of simplicial homology *****
ZZ + ZZ^2
> add monomial x^2 x*y y*z
> char 0
New characteristic: 0
> denom
1 - x^2*ZZ^2 - x*y*ZZ^2 - y*z*ZZ^2 - x^2*y*ZZ^3 - x*y*z*ZZ^3
> set multigrade false
> denom
1 - 3*ZZ^2 - 2*ZZ^3
> quit

Thanks for visiting.
\end{verbatim}
\caption{Calculation with {\tt poincare}}\label{ex:poincare}
\end{example}

Among the things we may observe in example \ref{ex:poincare} is the basic set of
commands -- {\tt add simplex} and {\tt add monomial} to build
simplicial complexes or monomial ideals, {\tt homology} and {\tt
denominator} (or an abbreviation {\tt denom} thereof) to calculate
simplicial homology and the Poincaré-Betti denominator respectively,
as well as the command {\tt char}, which changes the field
characteristic over which all homology calculations take place and the
command sequence {\tt set multigrade false}, which sets a flag that
causes the program to change the way it prints the polynomials output
by the {\tt denominator} command, so that instead of the polynomial
$b_{R,k}(\mathbb x, t)$ the program prints the polynomial
$b_{R,k}(1,\dots,1,t)$. The program ends upon receiving {\tt quit}.

In addition to these, there are the commands {\tt clear}: clearing the stored
simplicial complex and monomial ideal, but not changing the
characteristics used and {\tt var}: which changes the implicit
homology variable, which in my review in this paper has been called
$t$, and which by default in {\tt poincare} is called {\tt ZZ}. Should
you wish to use {\tt ZZ} as a ring variable in your calculations, a
change of homology variable will be necessary. For this, the command
{\tt var} is provided, with which you can change the string that 
{\tt poincare} uses for the homology variable.

\section{Questions and future directions}

There are several things that I want to improve upon on the system
herein presented, and also several questions that can be posed.

There are numerous complexity issues associated to the program in its
current form. Mainly, these issues are related to the size of the
resulting po-graphs for larger sets of monomials. As an example, the
initial ideal of a Gröbner basis with revlex ordering of the
homogenized cyclic 6-root ideal, i.e. the ideal generated by
\begin{align*}
&x_1+x_2+x_3+x_4+x_5+x_6,\\
&x_1x_2+x_2x_3+x_3x_4+x_4x_5+x_5x_6+x_1x_6,\\
&x_1x_2x_3+x_2x_3x_4+x_3x_4x_5+x_4x_5x_6+x_1x_5x_6+x_1x_2x_6,\\
&x_1x_2x_3x_4+x_2x_3x_4x_5+x_3x_4x_5x_6+
 x_1x_4x_5x_6+x_1x_2x_5x_6+x_1x_2x_3x_6,\\
&x_1x_2x_3x_4x_5+x_2x_3x_4x_5x_6+x_1x_3x_4x_5x_6+\\
&\quad x_1x_2x_4x_5x_6+x_1x_2x_3x_5x_6+x_1x_2x_3x_4x_6,\\
&x_1x_2x_3x_4x_5x_6-y^6
\end{align*}
invariably becomes larger than Pari's working memory. This initial ideal
has 100 monomial generators and produces a po-graph with 11443
elements. The calculations normally halt after between 200 and 400
lattice point calculations.

The problem I have observed with for instance this example is that the
calculation of homology of large simplicial complexes is memorywise unfeasible.
The most visible
problem is when the homology calculations turn out to be too hard,
since this results in a crash in the Pari library; whereas too hard
construction of the simplicial complexes merely result in slow running
of the program.

\begin{question}
Given the rather special structure of the matrices that are used to
calculate field homology of a simplicial complex, can anything be said
about the sizes of elements of a matrix as some sort of reduction
algorithm is used to deduce its rank? Can some variant of
row-reduction be found such that the matrices that produce ``large''
entries (for instance, larger than a standard 32-bit word) have some
easily recognizable feature?
\end{question}

If such a method could be found, then it would be possible to single
out the specific homology matrices that would actually need treatment
with some sort of bignum library, and produce fast 32-bit (sparse)
matrix arithmetic to deal with all other instances. Such a separate
treatment would also obliterate the need to rely on Pari's internal
memory allocation structures and would give the user more control over
acceptable memory consumption for the calculations at hand.

An obvious further point of attack is the generation algorithm for the
Berglund complexes. This would be vastly improved if a search
algorithm would be constructed that minimizes the number of multiple
checks done on each candidate monomial; since the graph connectivity
checks are not, in the context, very fast.

Finally, an interesting direction to take would be to look at the APIs
for larger computer algebra systems and try to adapt the code here
written to work as a pluggable module to those systems; for instance
providing an interface to calculate Poincaré-Betti series of monomial
rings directly from Singular or Magma or Macaulay 2.

\bibliographystyle{apalike}
\bibliography{../library}

\end{document}